\documentclass[10pt, article, titlepage,openany, rotate]{amsart}

\usepackage{amssymb, amsfonts,amscd,verbatim}
\usepackage{latexsym}

\usepackage{graphicx}
\usepackage{verbatim}   % useful for program listings
\usepackage{color}      % use if color is used in text
\usepackage{subfigure}  % use for side-by-side figures

\usepackage{amsmath}    % need for subequations

\theoremstyle{plain}
\newtheorem{Prop}{Proposition}[section]
\newtheorem{Thm}[Prop]{Theorem}
\newtheorem{Cor}[Prop]{Corollary}

\newtheorem{Lem}[Prop]{Lemma}

\theoremstyle{definition}
\newtheorem{Def}[Prop]{Definition}

\theoremstyle{remark}

\newtheorem{Problem}[Prop]{\bf Problem}

\newcommand{\map}{\operatorname{Map}\nolimits}

\def\dim{\mathop{\roman{dim}}}
\def\int{\mathop{\roman{int}}}

\def\1{^{-1}}
\def\Tilde{\widetilde}

\def\dim{\text{dim}}

\def\Tor{\text{Tor}}

\def\Ab{\text{Ab}}

\def\DDD{{\mathcal D}}
\def\PPP{{\mathcal P}}
\def\NNN{{\mathcal N}}
\def\AP{{\mathcal AP}}

\def\dokaz{{\bf Proof. }}
\numberwithin{equation}{section}
\newcommand{\cc}[1]{\!\!\!<\!\!#1\!\!>}
\begin{document}
%\fontsize{20}{34pt}\selectfont
\title[
Algebraic properties of quasi-finite complexes
]%
   {Algebraic properties of quasi-finite complexes}

\author{M.~Cencelj}
\address{Fakulteta za Matematiko in Fiziko,
Univerza v Ljubljani,
Jadranska ulica 19,
SI-1111 Ljubljana,
Slovenija }
\email{matija.cencelj@guest.arnes.si}

\author{J.~Dydak}
\address{University of Tennessee, Knoxville, TN 37996, USA}
\email{dydak@math.utk.edu}

\author{J.~Smrekar}
\address{Departament d'Algebra i Geometria,
Universitat de Barcelona,
Gran Via de les Corts Catalanes, 585,
E-08007 Barcelona,
Espa\~na}
\email{smrekar@ub.edu}

\author{A.~Vavpeti\v c}
\address{Fakulteta za Matematiko in Fiziko,
Univerza v Ljubljani,
Jadranska ulica 19,
SI-1111 Ljubljana,
Slovenija }

\email{ales.vavpetic@fmf.uni-lj.si}

\author{\v Z.~Virk}
\address{Fakulteta za Matematiko in Fiziko,
Univerza v Ljubljani,
Jadranska ulica 19,
SI-1111 Ljubljana,
Slovenija }

\email{ziga.virk@student.fmf.uni-lj.si}

\date{ September 23, 2005}


\thanks{Supported in part by the Slovenian-USA research grant BI--US/05-06/002 and the ARRS
research project No. J1--6128--0101--04}

\begin{abstract}

A countable CW complex $K$ is quasi-finite (as defined by
A.Karasev \cite{K}) if for every finite subcomplex $M$ of $K$
there is a finite subcomplex $e(M)$ such that any
map $f:A\to M$, where $A$ is closed in a separable metric space $X$
satisfying $X\tau K$, has an extension $g:X\to e(M)$.
Levin's \cite{Lev} results imply that none of the
Eilenberg-MacLane spaces $K(G,2)$ is quasi-finite if $G\ne 0$. In
this paper we discuss quasi-finiteness of all Eilenberg-MacLane
spaces. More generally, we deal with CW complexes with finitely many nonzero Postnikov invariants.

Here are the main results of the paper:
\begin{Thm}
Suppose $K$ is a countable CW complex
  with finitely many nonzero Postnikov invariants. If $\pi_1(K)$ is a locally finite group and $K$ is quasi-finite, then
$K$ is acyclic.
\end{Thm}

\begin{Thm}
Suppose $K$ is a countable non-contractible
CW complex
  with finitely many nonzero Postnikov invariants. If $\pi_1(K)$ is nilpotent and $K$ is quasi-finite, then
$K$ is extensionally equivalent to $S^1$.
\end{Thm}

\end{abstract}

\maketitle

\medskip
%Printed on \today.
\medskip
\tableofcontents

\section{Introduction}

\par The notation $K\in AE(X)$ or $X\tau K$ means that
any map $f:A\to K$, $A$ closed in $X$, extends over $X$.

\begin{Thm}[Chigogidze] \label{ChigEquivThm}
For each countable simplicial complex
$P$ the following conditions are equivalent:
\begin{itemize}
\item [1.]  $P\in AE(X)$ implies $P\in AE(\beta(X))$ for any normal
space $X$.

\item[2.]  There exists a $P$-invertible map $p:X\to I^\omega$ of a metrizable compactum
$X$ with $P\in AE(X)$ onto the Hilbert cube.
\end{itemize}
\end{Thm}

Karasev  \cite{K} gave an intrinsic characterization of countable complexes $P$ satisfying \ref{ChigEquivThm} and
called them {\it quasi-finite complexes}.
\begin{Def}\label{QuasiFiniteDef}
A CW complex $K$ is called {\it quasi-finite} if there is a function $e$ from the family
of all finite subcomplexes of $K$ to itself satisfying the following property:
For every separable metric space $X$ such that $K$ is an absolute extensor of $X$
 and for every map
$f:A\to M$, $A$ closed in $X$,  $f$ extends to $g:X\to e(M)$.\end{Def}

For subsequent generalizations of quasi-finiteness see  \cite{KV1}
and \cite{C-V}. In particular, it is shown in \cite{C-V} that
a countable CW complex $K$ is quasi-finite if and only if $X\tau K$
implies $\beta(X)\tau K$ for all separable metric spaces $X$. That is an improvement
of \ref{ChigEquivThm}.
\par The first example of a non-quasi-finite CW complex was given
by Dranishnikov \cite{DR2} who showed that $K(Z,4)$ admits a
separable metric space $X$ satisfying $X\tau K(Z,4)$ but not
$\beta(X)\tau K(Z,4)$ (see \cite{DW} for other examples of such
$X$). In \cite{Dy5} it was shown that all $K(G,n)$, $n\ge 3$ and
$G\ne 0$, admit a separable metric space $X$ so that $\dim_G(X)=n$
but $\dim_G(\beta(X)) > n$ (see also \cite{Kar} for related
resuts). Finally, Levin \cite{Lev} established a result implying
the same fact for all $K(G,n)$ so that $G\ne 0$ and $n\ge 2$. The
only remaining case among Eilenberg-Maclane spaces are complexes
$K(G,1)$.

\begin{Problem}\label{AreKGOneQFQ}  Characterize groups $G$ such that $K(G,1)$ is quasi-finite.
What are the properties of the class of groups $G$ such that
$K(G,1)$ is quasi-finite?
\end{Problem}

Problem \ref{AreKGOneQFQ} was the main motivation of this paper.
More generally, we discuss quasi-finiteness of complexes with
finitely many non-trivial Postnikov invariants.

\section{Truncated cohomology}

One of the main tools of this paper is truncated cohomology used
for the first time by Dydak and Walsh \cite{DW2} in their
construction of an infinite-dimensional compactum $X$ of integral
dimension $2$.
\par

Given a pointed CW complex $L$ and a pointed space $X$ we define
$h_L^k(X)$ as the $(-k)$-th homotopy group of the function space
$\map_*(X,L)$, the space of base-point preserving maps whose base-point is the constant map. Since we are interested in Abelian groups, $k$
ranges from minus infinity to $-2$. Also, spaces $X$ of interest
in this paper are countable CW complexes.

\par CW complexes $L$ for which trunctated cohomology $h^\ast_L$
is of most use are those with finite homotopy groups. In that case
$h^\ast_L$ is {\it continuous} in the sense that any map $f:K\to
\Omega^k L$ that is phantom (that means all restrictions $f|_M$
are homotopically trivial for finite subcomplexes $M$ of $K$) must
be homotopically trivial if $K$ is a countable CW complex. In case of $L$ having finite homotopy
groups, Levin \cite{LevFirst} (see Proposition 2.1) proved that
$h^\ast_L$ is {\it strongly continuous}: any map $f:N\to \Omega^k
L$, $N$ being a subcomplex of $K$, that cannot be extended over a
countable CW complex $K$, admits a finite subcomplex $M$ of $K$
such that $f|_{M\cap N}$ cannot be extended over $M$.
\par Since we are interested in vanishing of truncated cohomology
$h^\ast_L$, the remainder of this section is devoted to weak
contractibility of mapping spaces.
\par

We first recall a result that in the literature is known as the
Zabrodsky Lemma (see Miller \cite{Miller}, Proposition 9.5, and
Bousfield \cite{Bousfield}, Theorem 4.6 as well as Corollary 4.8).
\begin{Lem} \label{zabrodsky_lemma}
Let $F\to E\to B$ be a fibration where $B$ has the homotopy type
of a connected CW complex. Let $X$ be a space. If $\map_*(F,X)$ is
weakly contractible, the induced map $\map_*(B,X)\to\map_*(E,X)$
is a weak homotopy equivalence. \hfill $\blacksquare$
\end{Lem}

\begin{Def}\label{MillerComplexDef}
Let ${\mathcal P}$ be a set of primes.
By a {\it ${\mathcal P}$-complex} we mean a finite CW complex $K$ that
is simply connected and all its homotopy groups are ${\mathcal P}$-groups. That is,
homotopy groups of $K$ are finite and the order of each element is a product
of primes belonging to ${\mathcal P}$.
\par
A CW complex $K$ is a {\it co-${\mathcal P}$-complex} if for some $k$ the mapping space
$\map_*\big(\Sigma^kK,L\big)$ is weakly contractible for all ${\mathcal P}$-complexes $L$.
\end{Def}

\begin{Lem} \label{miller_zabrodsky_sullivan_mcgibbon_friedlander_mislin}
 If $K$ is one of the
following
\begin{enumerate}
\item   The classifying space $BG$ of a Lie group $G$ with a
finite number of path components,

\item   A connected infinite loop space whose fundamental group is
a torsion group,

\item   A simply connected CW complex with finitely many homotopy
groups,
\end{enumerate}
then $\map_*(K,L)$ is weakly contractible for all nilpotent finite
complexes $L$ with finite homotopy groups.
\end{Lem}
\dokaz Let $L$ be a finite nilpotent complex with finite homotopy groups. 
The hypotheses render $L$
complete with respect to Sullivan's finite completion (see
\cite{Sullivan}). Thus case (1) follows from Friedlander and
Mislin \cite{Friedlander-Mislin}, Theorem 3.1, while case (2)
follows from McGibbon \cite{McGibbon}, Theorem 3. Case (3) follows
from \ref{zabrodsky_lemma} and (2) by induction over the number of
non-trivial homotopy groups of $K$. See more details in the proof
of \ref{miller_zabrodsky_sullivan}.
 \hfill
$\blacksquare$

\begin{Prop}\label{FiniteWedgeOfcoPcomplexes}
A finite product (or a finite wedge) of co-$\PPP$-complexes is a
co-$\PPP$-complex.
\end{Prop}
\dokaz In case of a finite wedge the proof is quite simple as
$\map_\ast(K\vee P,L)$ is the product of $\map_\ast(K,L)$ and $\map_\ast(P,L)$.
For the finite product one can use induction plus an observation
that \ref{zabrodsky_lemma} can be applied to a fibration $F\to
K\to P$ and yield that $K$ is a co-$\PPP$-complex if both $F$ and
$P$ are co-$\PPP$-complexes.
 \hfill $\blacksquare$

 \begin{Prop}\label{ArbitraryWedgeOfcoPcomplexes}
 Let ${\mathcal P}$ be a set of primes.
 Suppose $K_s$, $s\in S$, is a a family of CW complexes. If
 there is a natural number $k$ so that all function spaces
 $\map_\ast(\Sigma^k K_s,L)$ are weakly contractible for all $\PPP$-complexes $L$, then the wedge
 $K=\bigvee\limits_{s\in S}K_s$ is a
co-$\PPP$-complex. Moreover, if $S$ is countable and each $K_s$ is
countable, then the weak product $\prod\limits_{s\in S}K_s$ is a
co-$\PPP$-complex.
\end{Prop}
\dokaz The case of the wedge is left to the reader. If $S$ is
countable, then each finite product $K_T=\prod\limits_{s\in T}K_s$
has the property that $\map_\ast(K_T,\Omega^kL)$ is weakly contractible
for any $\PPP$-complex $L$ as in the proof of
\ref{FiniteWedgeOfcoPcomplexes}. Using the fact that truncated
cohomology with respect to $\Omega^kL$ is continuous, one gets
that $K'=\prod\limits_{s\in S}K_s$, being the direct limit of
$K_T$, also has the property that $\map_\ast(K',\Omega^kL)$ is weakly
contractible.
 \hfill $\blacksquare$

\begin{Def}\label{WeakMillerGroupDef}
Let ${\mathcal P}$ be a set of primes and let $G$ be a group. $G$
is called a {\it co-${\mathcal P}$-group} if $K(G,1)$ is a
co-${\mathcal P}$-complex.
\end{Def}

% The name stems from the celebrated resolution of the Sullivan conjecture by Miller (see \cite{Miller}).

 By
Miller's theorem, all locally finite groups are co-${\mathcal P}$-groups,
where ${\mathcal P}$ is the set of all primes. Another example would consist of all acyclic
groups. Divisible groups would serve as well. Note that by the Zabrodsky Lemma \ref{zabrodsky_lemma}, a group extension
$N\rightarrowtail G\twoheadrightarrow Q$ implies that under the assumption that $N$ is a co-${\mathcal P}$-group,
$G$ is a co-${\mathcal P}$-group if and only if $Q$ is.

\begin{Def}
Let $K$ be a connected CW complex. We say that $K$ has {\it
finitely many unstable Postnikov invariants} if for some $k\ge 0$,
the $k$-connected cover $K\cc{k}$ of $K$ is an infinite loop
space. As usual, $K\cc{k}$ is the (homotopy) fibre of the $k$-th
Postnikov approximation $K\to P_k(K)$.
\end{Def}

Note that infinite loop spaces (in particular infinite symmetric
products) and Postnikov pieces are special cases.

\begin{Lem} \label{miller_zabrodsky_sullivan}
Suppose ${\mathcal P}$ is a set of primes. Let $K$ be a connected
CW complex with finitely many unstable Postnikov invariants. $K$
is a co-${\mathcal P}$-complex if and only if $G=\pi_1(K)$ is a
co-${\mathcal P}$-group.
\end{Lem}
\dokaz Let $L$ be a $\PPP$-complex. Let $\Tilde K$ be the
universal cover of $K$. If $K$ is itself an infinite loop space,
so is $\Tilde K$, and therefore the space $\map_*(\Tilde K,L)$ is
weakly contractible by Theorem 3 of McGibbon \cite{McGibbon}.
Otherwise for some $i\ge 1$ the $i$-connected cover $\Tilde
K\cc{i}$ of $\Tilde K$ is an infinite loop space. Consider the
fibration sequence $\Tilde K\cc{i}\to\Tilde K\to P_i\Tilde K$
where $P_i\Tilde K$ is the $i$-th Postnikov approximation of
$\Tilde K$. The space $\map_*(\Tilde K\cc{i},L)$ is weakly
contractible by Theorem 3 of McGibbon \cite{McGibbon}. It follows
essentialy from Zabrodsky \cite{Zabrodsky}, Theorem D, and the
fact that $L$ is Sullivan-complete, that the mapping space
$\map_*(P_i\Tilde K,L)$ is weakly contractible (see also McGibbon
\cite{McGibbon}, Theorem 2). Thus by Lemma \ref{zabrodsky_lemma},
also the space $\map_*(\Tilde K,L)$ is weakly contractible. The
space $\Tilde K$ sits in the fibration sequence $\Tilde K\to K\to
K(G,1)$ and another application of Lemma \ref{zabrodsky_lemma}
renders the spaces $\map_*(K(G,1),L)$ and $\map_*(K,L)$ weakly
equivalent.
 \hfill $\blacksquare$

\begin{Lem} \label{local_away_from_P_to_P_torsion}
Let $\PPP$ be a nonempty set of primes. If $G$ is a nilpotent
group that is local away from $\PPP$, then it is a
co-$\PPP$-group.
\end{Lem}
\dokaz Let $\PPP'$ denote the set of primes not in $\PPP$. The
hypotheses on $G$ render $K(G,1)$ a $\PPP'$-local space. By the
fundamental theorem of localization of nilpotent spaces it follows
that the homology of $K(G,1)$ is also $\PPP'$-local. Let $\dots\to
L_3\to L_2\to L_1\to L_0$ denote the refined Postnikov tower for
$L$. That is, $L_0$ is a point and for each $i$, the fibration
$L_i\to L_{i-1}$ is principal with fibre $K(G_i,k_i)$ where $G_i$
is $p$-torsion abelian. Note that $L$ is weakly equivalent to the
inverse limit $\lim_iL_i$, and since $K(G,1)$ is a CW complex it
suffices to show that $\map_*\big(K(G,1),\lim_iL_i\big)$ is weakly
contractible. This latter space is homeomorphic with the inverse
limit $\lim_i\map_*(K(G,1),L_i)$. Since the fibrations are
principal, the Puppe sequence shows that we only need to consider
reduced cohomology $\Tilde H^*(K(G,1);G_i)$ with coefficients in
$G_i$. Since $H_*(G)$ is local away from $\PPP$ it follows by the
universal coefficient theorem that $\Tilde H^*(K(G,1);G_i)$ is
trivial. \hfill $\blacksquare$

\begin{Cor} \label{NilpotentCo-P-groups}
Suppose ${\mathcal P}$ is a set of primes and $G$ is a nilpotent
group with Abelianization $\Ab(G)$.
 If $\Ab(G)/\Tor(\Ab(G))$ is ${\mathcal P}$-divisible, then $G$ is a co-${\mathcal P}$-group.
\end{Cor}
\dokaz By \ref{bridge} $\Ab(G)/\Tor(\Ab(G))$ is ${\mathcal
P}$-divisible if and only if $G$ is local away from $\PPP$. \hfill
$\blacksquare$

\section{Homology and cohomology of quasi-finite CW complexes}

In this section we deal with (co)homological properties of
quasi-finite complexes. First, we need a generalization of Theorem
II of \cite{Dy6}.

\begin{Thm} \label{HomologyRealizationThm}
Suppose $K$ is a countable CW complex
and $h_\ast$ is a generalized reduced homology theory such that
$h_\ast(K)=0$. For any CW complex $P$ and any $\alpha\in h_\ast(P)\setminus \{0\}$
there is a compactum $X$ and a map $f:A\to P$ from a closed subset $A$ of $X$
such that $X\tau K$, $\alpha=f_\ast(\gamma)$ for some $\gamma\in \check h_\ast(A)$
and $\gamma$ is $0$ in $\check h_\ast(X)$.
\end{Thm}
\dokaz Replacing $P$ by the carrier of $\alpha$ we may assume $P$
is finite. Compactum $X$ is built as in Theorem II of \cite{Dy6}.
We start with $X_1=Cone(P)$, $A_1=P$ and build an inverse sequence
$(X_n,A_n)$ of compact polyhedra so that for every extension
problem $g:B\to K$, $B$ closed in $X_n$, there is $m > n$ and a
map $G:X_m\to K$ extending $g\circ p^m_n:B'\to K$, where
$p^m_n:X_m\to X_n$ is the bonding map and $B'=(p^m_n)^{-1}(B)$.
For each $n$ we have $\gamma_n\in h_\ast(A_n)$ which vanishes in
$h_\ast(X_n)$. In the inductive step we pick an extension problem
$g:B\to K$, $B$ closed in $X_n$, create an extension $G:X_n\to
Cone(K)$, and consider the pull-back $E$ of the projection
$K\times I\to Cone(K)$ under $G$. The projection $p:E\to X_n$ has
fibers being either homeomorphic to $K$ or single points.
Therefore $h_\ast(p)$ is an isomorphism and one can pick a finite
subpolyhedron $A_{n+1}$ of $E$ carrying $\gamma_{n+1}\in
h_\ast(A_{n+1})$ which gets mapped to $\gamma_n$ under
$h_\ast(p)$. Since $\gamma_{n+1}$ vanishes in $h_\ast(E)$, it
vanishes in a finite subpolyhedron $X_{n+1}$ of $E$ containing
$A_{n+1}$. Since there are only countably many extension problems
to be solved (see \cite{DR1} or \cite{DR4}) that process produces
an inverse sequence whose inverse limit $(X,A)$ satisfies $X\tau
K$ and one has $\gamma\in\check h_\ast(A)$ so that $\gamma$
vanishes in $\check h_\ast(X)$ and $f_\ast(\gamma)=\alpha$, where
$f:A\to P=A_1$ is the projection. \hfill $\blacksquare$

\begin{Thm} \label{CohomologyRealizationTheorem}
Suppose $K$ is a countable CW complex and $h^\ast$ is a
strongly continuous truncated cohomology theory such that $h^\ast(K)=0$. For
any countable CW complex $P$ and any $\alpha\in h^\ast(P)\setminus
\{0\}$ there is a compactum $X$ and a map $f:A\to P$ from a closed
subset $A$ of $X$ such that $X\tau K$ and there is no $\gamma\in
\check h^\ast(X)$ satisfying $\gamma|_A=f^\ast(\alpha)$.
\end{Thm}
\dokaz We can reduce the proof to the case of $P$ being a finite
polyhedron as there is a finite subcomplex $M$ of $P$ so that
$\alpha|_M\ne 0$ and that $M$ can be used instead of $P$.
Compactum $X$ is built as in \ref{HomologyRealizationThm}. We
start with $X_1=Cone(P)$, $A_1=P$ and built an inverse sequence
$(X_n,A_n)$ of compact polyhedra so that for every extension
problem $g:B\to K$, $B$ closed in $X_n$, there is $m > n$ and a
map $G:X_m\to K$ extending $g\circ p^m_n:B'\to K$, where
$p^m_n:X_m\to X_n$ is the bonding map and $B'=(p^m_n)^{-1}(B)$.
Also, for each $n$ the pullback $\alpha_n$ of $\alpha$ under
$A_n\to A_1$ does not extend over $X_n$. In the inductive step we
pick an extension problem $g:B\to K$, $B$ closed in $X_n$, create
an extension $G:X_n\to Cone(K)$, and consider the pull-back $E$ of
the projection $K\times I\to Cone(K)$ under $G$. The projection
$p:E\to X_n$ has fibers being either homeomorphic to $K$ or single
points. Therefore $p^\ast=h^\ast(p)$ is an isomorphism. Since
$p^\ast(\alpha_n)$ does not extend over $E$, there is a finite
subpolyhedron $X_{n+1}$ of $E$ such that $p^\ast(\alpha_n)$
restricted to $A_{n+1}=X_{n+1}\cap p^{-1}(A_n)$ does not extend
over $X_{n+1}$. Since there are only countably many extension
problems to be solved (see \cite{DR1} or \cite{DR4}) that process
produces an inverse sequence whose inverse limit $(X,A)$ satisfies
$X\tau K$ and the projection $f:A\to P=A_1$ has the property that
there is no $\gamma\in \check h^\ast(X)$ satisfying
$\gamma|_A=f^\ast(\alpha)$.
 \hfill $\blacksquare$

Recall that, given a map $i:M\to N$, $X\tau i$ means that for any
map $f:A\to M$, $A$ closed in $X$, there is a map $g:X\to N$
extending $i\circ f$.
\begin{Thm} \label{HomCohOfQF}
Suppose $K$ is a countable CW complex
and $i:M\to N$ is a map of CW complexes such that $X\tau K$ implies
$X\tau i$ for all compacta $X$.
\begin{enumerate}
\item If $h_\ast$ is a generalized reduced homology theory such that
 the inclusion induced homomorphism
$h_\ast(M)\to h_\ast(N)$ is not trivial, then
$h_\ast(K)\ne 0$.

\item If $h^\ast$ is a truncated strongly continuous cohomology theory such that
the inclusion induced homomorphism
$h^\ast(N)\to h^\ast(M)$ is not trivial  and $M$ is countable, then $h^\ast(K)\ne 0$.
\end{enumerate}
\end{Thm}
\dokaz We may assume $i$ is an inclusion.
\par
1. Suppose $\alpha\in h_\ast(M)$ does not become $0$ in $h_\ast(N)$.
As in \ref{HomologyRealizationThm} pick a map $f:A\to M$ of a closed
subset of a compactum $X$ so that $X\tau K$ and $\gamma$ equals $0$ in
$\check h_\ast(e(M))$ for some $\gamma\in \check h_\ast(A)$ satisfying
$f_\ast(\gamma)=\alpha$. If $f$ extends to $g:X\to N$,
then $\alpha=f_\ast(\gamma)$ becomes $0$ in $h_\ast(N)$, a contradiction.
\par 2. Suppose $\alpha\in h^\ast(N)$ and $\alpha|_M\ne 0$.
We may reduce this case to $M$ finite by switching to a finite
subcomplex $L$ of $M$ with the property $\alpha|_L\ne 0$.
As in \ref{CohomologyRealizationTheorem} pick a map $f:A\to M$
of a closed subset of a compactum $X$ so that $X\tau K$ and
$f^\ast(\alpha|_M)$ does not extend over $X$.
If $f:A\to M$ extends to $g:X\to N$, then $g^\ast(\alpha)\in \check h_\ast(X))$
extends $f^\ast(\alpha|_M)$, a contradiction.
\hfill $\blacksquare$

\begin{Thm}\label{coP-QFAreAcyclic}
Suppose ${\mathcal P}$ is a set of primes.
Let $K$ be a connected countable co-${\mathcal P}$-complex.
If $K$ is quasi-finite, then it is $Z_{({\mathcal P})}$-acyclic.
\end{Thm}
\dokaz Assume $K$ is quasi-finite and not $Z_{({\mathcal
P})}$-acyclic. Replace $K$ with $\Sigma K$ (using \cite{C-V}) if necessary to ensure
$H_k(K;Z_{({\mathcal P})})\neq 0$ for some $k\ge 2$. Let
$\alpha_K\in H_k(K;Z_{({\mathcal P})})$ be nonzero. Since $K$ is
the colimit of its finite subcomplexes, $\alpha_K$ is the image of
$\alpha_M\in H_k(M;Z_{({\mathcal P})})$ for some finite subcomplex
$M$ of $K$. Certainly the image of $\alpha_M$ under
$H_k(M;Z_{({\mathcal P})})\to H_k(e(M);Z_{({\mathcal P})})$ is
nontrivial. Thus Lemma  \ref{LevinStyleHomologyLemma} yields a
${\mathcal P}$-complex $L$ with the restriction morphism
$[e(M),\Omega^2L]\to[M,\Omega^2L]$ nontrivial. This is to say that
$h^*(e(M))\to h^*(M)$ is nontrivial where $h^*$ is the truncated
cohomology theory defined by virtue of $\Omega^2L$. The hypotheses
on $L$ ensure strong continuity of $h^*$. Thus the nontriviality
of $h^*(e(M))\to h^*(M)$ contradicts (2) of Theorem
\ref{HomCohOfQF}. \hfill $\blacksquare$

\begin{Cor}\label{pi_1_K_Not_Divisible_If_K_Postnikov_is_QF}
Let $K$ be a countable CW complex with finitely many nonzero
homotopy groups and $G=\pi_1(K)$ nilpotent. Suppose that $G$ is
not torsion. If $K$ is quasi-finite, the group $FG=G/\Tor(G)$ (and
thus also $\Ab(G)/\Tor(\Ab(G))$) is not divisible by any prime
$p$.
\end{Cor}
\dokaz Suppose that, on the contrary, $FG$ is divisible by a prime
$p$, hence local away from $p$. Since $G$ is not torsion and is
nilpotent, also $\Ab(G)$ is not torsion, hence certainly
$H_1(K)\otimes Z_{(p)}$ is nontrivial. Thus Theorem
\ref{coP-QFAreAcyclic} yields a contradiction. \hfill
$\blacksquare$

\ref{coP-QFAreAcyclic} and \ref{miller_zabrodsky_sullivan} imply the following.

\begin{Cor}
Let $K$ be a simply connected countable CW complex with at least one and at most finitely many nontrivial
homotopy groups. Then $K$ is not quasi-finite. \hfill $\blacksquare$
\end{Cor}

\begin{Cor}\label{LocFiniteKGOne}
Suppose $G$ is a locally finite countable group. If $K(G,1)$ is quasi-finite, then $G$ is acyclic. \hfill $\blacksquare$
\end{Cor}

However, there are some countable acyclic groups $G$ for which $K(G,1)$ are
also not quasi-finite. Cencelj and Repov\v s \cite{CeRe}, using results of
Dranishnikov and Repov\v s \cite{DrR1} showed in \S 5 that the minimal
grope $M^*$ which is $K(\pi_1(M^*),1)$ is not quasi-finite.
This holds also for the fundamental
group of any grope: For a grope $M$ let $\gamma(m)$ denote the maximal number of handles
on the discs with handles used in the construction of the $m$-stage of $M$. Modify the
inverse limit construction of the example of \cite{CeRe} replacing every simplex in the triangulation
of the $k$-th element of the inverse system by the $n$-th stage of the grope which has
every generator replaced by a disc with $\gamma(kn)$ handles.

\section{Ljubljana complexes}

\begin{Def}\label{LjubljanaComplexDef}
A connected CW complex $L$ is called a {\it Ljubljana complex} if there is
a co-$\AP$-complex $K$, $\AP$ being the set of all primes, such that, for any compactum $X$,
the conditions $X\tau L$ and $X\tau K(H_1(K),1)$ imply $X\tau K$.
\end{Def}

\begin{Lem} \label{FibrationOfLjuComplexes}
Suppose $F\to E\to B$ is a fibration of connected CW complexes.
If $F$ is a co-$\AP$-complex, $\AP$ being the set of all primes, and $B$ is a Ljubljana complex, then $E$ is a Ljubljana complex.
\end{Lem}
\dokaz Notice $\pi_1(E)\to \pi_1(B)$ is an epimorphism (use the long exact sequence of a fibration) which implies $H_1(E)\to H_1(B)$ is an epimorphism.

Pick a co-$\AP$-complex $K$ such that $X\tau K$
and $X\tau K(H_1(B),1)$ imply $X\tau B$ for all compacta $X$.
Let $M$ be the wedge of $F$, $K$, $K(Q,1)$, and $X\tau K(Z/p^\infty,1)$ for all primes $p$. By \ref{ArbitraryWedgeOfcoPcomplexes} and the Miller Theorem, $M$ is a co-$\AP$-complex.
Suppose $X$ is a compactum such that $X\tau M$ and $X\tau K(H_1(E),1)$.
By \ref{EpiAndCohDim} one gets $X\tau K(H_1(B),1)$ which, together with $X\tau K$,
implies $X\tau B$. Since $X\tau F$ and $X\tau B$, we infer $X\tau E$.
\hfill $\blacksquare$

\begin{Cor} \label{ExamplesOfCComplexes}
Let $L$ be a connected CW complex with nilpotent fundamental
group. If $L$ has finitely many unstable Postnikov invariants,
 then $L$ is a Ljubljana complex.
\end{Cor}
\dokaz 
Notice that the universal cover $\Tilde L$ of $L$ is a co-$\AP$-complex by \ref{miller_zabrodsky_sullivan}. We get $K(\pi_1(L),1)$ is a Ljubljana complex by \ref{CohDimComparisonGvsAbG}.
The fibration $\Tilde L\to L\to K(\pi_1(L),1)$ implies
$L$ is a Ljubljana complex.
\hfill $\blacksquare$

\begin{Def}\label{eAbelianComplexDef}
A connected CW complex $L$ is called {\it extensionally Abelian}
if $X\tau K(H_n(L),n)$ for all $n\ge 1$ imply $X\tau K$ for all
compacta $X$.
\end{Def}

\begin{Prop}\label{eAbelianAreLjuComplexes}
Each extensionally Abelian complex $L$ is a Ljubljana complex.
\end{Prop}
\dokaz Let $K$ be the weak product of $K(H_n(L),n)$, $n\ge 2$.
By (2) of \ref{miller_zabrodsky_sullivan_mcgibbon_friedlander_mislin}, $K$ is a co-$\AP$-complex. Clearly, $X\tau K$ and $X\tau K(H_1(L),1)$
imply $X\tau K(H_n(L),n)$ for all $n\ge 1$. Thus $X\tau L$.
\hfill $\blacksquare$

\begin{Prop}\label{WedgeOfLjuComplexes}
A finite wedge (or finite product) of Ljubljana complexes is a Ljubljana
complex.
\end{Prop}
\dokaz Let $L$ be the wedge (or the product) of Ljubljana
complexes $L_s$, $s\in S$, where $S$ is finite. For each $s\in S$ choose a co-$\AP$-complex $K_s$
such that for any compactum $X$ the conditions $X\tau K_s$ and $X\tau K(H_1(L_s),1)$
imply $X\tau L_s$. Let $K$ be the wedge of all $K_s$. By \ref{FiniteWedgeOfcoPcomplexes} it is a co-$\AP$-complex.
Notice that $H_1(L_s)$ is a retract of $H_1(L)$ 
for each $s\in S$. Therefore any compactum $X$
satisfying
\begin{itemize}
\item[a.] $X\tau K(H_1(L),1)$,

\item[b.] $X\tau K$,
\end{itemize}
also
satisfies $X\tau K(H_1(L_s),1)$ for each $s\in S$. Hence $X\tau
L_s$ for each $s\in S$ which implies $X\tau L$. \hfill
$\blacksquare$

There is a connection between Ljubljana complexes and
co-$\PPP$-complexes.

\begin{Prop}\label{LjuArecoP}
Suppose $K$ is a countable Ljubljana complex. If $\PPP$ is a set
of primes such that $H_1(K)/\Tor(H_1(K))$ is $\PPP$-divisible,
then $K$ is a co-$\PPP$-complex.
\end{Prop}
\dokaz Choose a co-$\AP$-complex $L$ such that, for any compactum $X$,
the conditions $X\tau L$ and $X\tau K(H_1(K),1)$ imply $X\tau K$.
Let $\PPP'$ be the complement of $\PPP$ in the set of all
primes. Consider $K'$, the wedge of $L$, $K(Z_{(\PPP')},1)$, $K(Q,1)$,
and all $K(Z/p,1)$ ($p$ ranging through all primes). By
\ref{NilpotentCo-P-groups} and \ref{ArbitraryWedgeOfcoPcomplexes} $K'$ is a co-$\PPP$-complex. Since
$X\tau K'$ implies $X\tau K$ for all compacta,  \ref{HomCohOfQF}
implies that there is $k\ge 0$ such that the truncated cohomology
of $K$ with respect to $\Omega^k L$, $L$ any $\PPP$-complex, is
trivial. Thus $K$ is a co-$\PPP$-complex. \hfill $\blacksquare$

\begin{Thm} \label{BigResultOnQF}
Suppose $K$ is a countable Ljubljana complex
such that $\Sigma^m K$ is equivalent (over the class of compacta) to a
quasi-finite countable complex $L$ for some $m\ge 0$. If $K$ is not acyclic, then it is
equivalent to $S^1$.
\end{Thm}
 \dokaz 
  We may assume $L$ is simply connected as $\Sigma^{m+1} K$
is equivalent to $\Sigma L$ (see \cite{DrDy2}) and $\Sigma L$ is
quasi-finite by \cite{C-V}.
\par Suppose $K$ is not equivalent to $S^1$. Choose a co-$\AP$-complex
$P$ such that conditions $X\tau P$ and $X\tau K(H_1(K),1)$ imply $X\tau K$.
Let $k\ge 2$ be a number such that all maps $\Sigma^n P\to R$
are null-homotopic if $R$ is an $\AP$-complex and $n\ge k$.
\par
 Step 1. $L$ is not
contractible as otherwise $\Sigma^m K$ would have to be
contractible implying $K$ being acyclic.
\par Step 2. $L$ is not acyclic as it is not contractible by Step 1.
\par Step 3.
Since $X\tau K$ implies $X\tau K(H_1(K),1)$, the group $H_1(K)$
has the property of $H_1(K)/\Tor(H_1(K))$ being divisible by some
prime $p$. Indeed, If $H_1(K)/\Tor(H_1(K))$ is not divisible by
any prime, then the Bockstein basis of $H_1(K)$ consists of all
Bockstein groups and $X\tau K(H_1(K),1)$ implies $X\tau S^1$ by
Bockstein First Theorem. Since $X\tau K$ implies $X\tau
K(H_1(K),1)$ and $X\tau S^1$ implies $X\tau K$ for any compactum
$X$, $K$ is equivalent to $S^1$ over compacta. 
\par Let $e$ be the
function of $L$.

\par Case 1: $H_\ast(K)$ is a torsion group. There is $M$ such that $H_\ast(M)\to H_\ast(e(M))$ is not trivial.
By \ref{LevinHomologyInFullStrength}, there is a map
$f:\Sigma^k(e(M))\to J$ such that $f|_{\Sigma^k(M)}$ is not
trivial, $J$ is simply connected, and all homotopy groups of $J$
are finite. Consider the wedge $N$ of $P$ and $K(\bigoplus\limits_{q} Z/q,1)$. 
Notice $X\tau N$ implies $X\tau
K(H_1(K),1)$. Therefore
 $X\tau N$ implies $X\tau K$ which, in turn, implies
$X\tau L$ and $X\tau i_M$, where $i_M:M\to e(M)$. Since
$\map_\ast(N,\Omega^k(J))$ is weakly contractible, \ref{HomCohOfQF}
implies homotopy triviality of $f|_{\Sigma^k(M)}$, a
contradiction.

\par Case 2: $H_\ast(K)$ is not a torsion group.
Notice $H_\ast(L)$ is not torsion as well. Indeed if $H_\ast(L)$
is torsion, we could find a finite dimensional compactum $Y$ of
high rational dimension but all torsion dimensions equal $1$. Such
compactum satisfies $Y\tau L$ but $Y\tau \Sigma^m K$ fails as it
implies rational dimension of $Y$ to be at most $m+n$, where
$H_n(K)$ is not torsion.
 There is $M$ such that
the image of $H_\ast(M)\to H_\ast(e(M))$ is not torsion. Therefore
there is $n> 0$  such that $H_n(M;Z_{(p)})\to H_n(e(M);Z_{(p)})$ is
not trivial. By \ref{LevinStyleHomologyLemma}, there is a map
$f:\Sigma^k(e(M))\to J$ such that $f|_{\Sigma^k(M)}$ is not
trivial, $J$ is simply connected, and all homotopy groups of $J$
are finite $p$-groups. Consider the wedge $N$ of $P$ and
$K(Z[\frac{1}{p}]\oplus Z/p,1)$. Using \ref{CohDimWRTNilpotentDp} and \ref{EpiAndCohDim}, one gets
 $X\tau N$ implies $X\tau K$ which, in turn, implies
$X\tau L$ and $X\tau i_M$, where $i_M:M\to e(M)$. Since
$\map_\ast(N,\Omega^k(J))$ is weakly contractible, \ref{HomCohOfQF}
implies homotopy triviality of $f|_{\Sigma^k(M)}$, a
contradiction. \hfill $\blacksquare$

\begin{Cor}\label{TorsionfreeKGOne}
Suppose $G$ is a nontrivial nilpotent group.
If $K(G,1)$ is quasi-finite, then it is equivalent, over the class
of paracompact spaces, to $S^1$.
\end{Cor}

\section{Application to cohomological dimension theory}

\begin{Thm} \label{UVW}
Suppose $G\ne 1$ is a countable group such that
$\dim_G(\beta(X))=1$ for every separable metric space $X$
satisfying $\dim_G(X)=1$. If $G$ is nilpotent, then $\dim_G(X)\leq
1$ implies $\dim(X)\leq 1$ for all paracompact spaces $X$.
\end{Thm}
\dokaz By an improvement of Chigogidze's Theorem \ref{ChigEquivThm}
contained in \cite{C-V}, $K(G,1)$ is quasi-finite. Therefore
\ref{TorsionfreeKGOne} says $K(G,1)$ is equivalent to $S^1$ over
compacta. A result in \cite{C-V} says that $K(G,1)$ is equivalent
to $S^1$ over paracompact spaces which completes the proof. \hfill
$\blacksquare$

\section{Appendix A}
In this section we discuss results related to groups that are needed in the paper.

\begin{Lem} \label{ImagesBelongToClass}
Let $p$ be a natural number and let ${\mathcal D}_p$ be the class of groups $G$
such that $\Ab(G)/Tor(\Ab(G))$ is $p$-divisible, where $\Ab(G)$ is the Abelianization of $G$.
 If $f:G\to H$ is an epimorphism and $G\in {\mathcal D}_p$, then $H\in {\mathcal D}_p$.
\end{Lem}
\dokaz Notice that $G\in \DDD_p$ if and only if for each $a\in G$ there is $b\in G$ and $k\ge 1$
such that $(a\cdot b^{-p})^k$ belongs to the commutator subgroup $[G,G]$ of $G$. Suppose $a\in H$. Pick $b\in G$ with $a=f(b)$. There is $c\in G$ and $k\ge 1$ such that
$(b\cdot c^{-p})^k\in [G,G]$. Now $(a\cdot f(c)^{-p})^k\in [H,H]$ and $H\in {\mathcal D}_p$.
\hfill $\blacksquare$

\begin{Lem} \label{TensorPrBelongsToClass}
Let $p$ be a natural number and let ${\mathcal D}_p$ be the class of groups $G$
such that $H/Tor(H)$ is $p$-divisible, where $H$ is the Abelianization of $G$.
 If $G,H$ are Abelian and $G\in {\mathcal D}_p$, then $G\otimes H\in {\mathcal D}_p$.
\end{Lem}
\dokaz It suffices to show
that for each element $a$ of $G\otimes H$ there is $b\in G\otimes H$
and an integer $k\ne 0$ such that $k\cdot a+kp\cdot b=0$.
That in turn can be reduced to generators of $G\otimes H$ of the form $g\otimes h$.
Pick $u\in G$ and an integer $k\ne 0$ such that $k\cdot g+kp\cdot u=0$.
Now $k\cdot (g\otimes h)+kp\cdot (u\otimes h)=0$.
\hfill $\blacksquare$

We recall a result of Robinson (see \cite{Robinson}, 5.2.6) on the
relation between a nilpotent group and its abelianization.

\begin{Prop}[Robinson] \label{robinson_lemma}
Let $\NNN$ denote the category of nilpotent groups. Let $\PPP$ be
a class of groups in $\NNN$ with the following properties.
\begin{enumerate}
\item   For $A$ and $B$ abelian, $B\in\PPP$, any quotient of
$A\otimes B$ belongs to $\PPP$. \item   For $K,Q\in\PPP$, an
extension $1\to K\to G\to Q\to 1$ in $\NNN$ implies $G\in\PPP$.
\end{enumerate}
Suppose that $G\in\NNN$. If $\Ab(G)$ belongs to $\PPP$, so does
$G$. \hfill $\blacksquare$
\end{Prop}

We note the following corollary.
\begin{Cor}\label{bridge}
Let $G$ be a nilpotent group and set $H=\Ab(G)$. If $H/\Tor(H)$ is
$p$-divisible, so is $G/\Tor(G)$.
\end{Cor}
\dokaz Define the class $\DDD_p$ by letting a nilpotent group $G$
belong to $\DDD_p$ if and only if $F_p(G)=G/\Tor_p(G)$ is
$p$-divisible where $\Tor_p(G)$ denotes the $p$-torsion subgroup
of $G$. Note that $F_p(G)$ is $p$-divisible if and only if
$G/\Tor(G)$ is, hence it suffices to check properties (1) and (2)
of Proposition \ref{robinson_lemma}.

As for (1) it follows from  \ref{ImagesBelongToClass} and \ref{TensorPrBelongsToClass}.

For (2), note that $F_p$ is a functor $\NNN\to\NNN$. Let $1\to
K\to G\to Q\to 1$ be an extension in $\NNN$. We apply $F_p$. Since
$\Tor_p(K)=K\cap\Tor_p(G)$, the morphism $F_p(K)\to F_p(G)$ is
injective. Evidently, $q\colon F_p(G)\to F_p(Q)$ is surjective.
Moreover, $F_p(K)$ is a subset of the kernel of $q$. Assume that
$K$ belongs to $\DDD_p$. If $q(\xi)=1$ for some $\xi\in F_p(G)$,
then $\xi^{p^i}\in F_p(K)$ for large enough $i$. By assumption on
$K$, the group $F_p(K)$ is $p$-divisible, hence
$\xi^{p^i}=\eta^{p^{i}}$ for an element $\eta\in F_p(K)$. But
$F_p(G)$ is free of $p$-torsion (and nilpotent), so the equality
$\xi^{p^i}=\eta^{p^{i}}$ in $F_p(G)$ implies $\xi=\eta$ (see for
example Hilton, Mislin, Roitberg \cite{HMR}, Corollary 2.3).
Therefore in fact $\xi\in F_p(K)$, ie $\ker q=F_p(K)$. This is to
say that $1\to F_p(K)\to F_p(G)\to F_p(Q)\to 1$ is an extension in
$\NNN$. If, in addition, $Q$ belongs to $\DDD_p$ then $F_p(K)$ and
$F_p(Q)$ are $p$-divisible and free of $p$-torsion, and as such
local away from $p$. Therefore so is $F_p(G)$, by Corollary 2.5 of
\cite{HMR}. \hfill $\blacksquare$

\begin{Lem} \label{EpiAndCohDim}
Let $f:G\to H$ be an epimorphism of Abelian groups and let $X$ be
a compactum. If

\begin{itemize}
\item[a.] $X\tau K(G,1)$,

\item[b.] $X\tau K(Q,1)$,

\item[c.] $X\tau K(Z/p^\infty,1)$ for all primes $p$,
\end{itemize}

then $X\tau K(H,1)$.
\end{Lem}
\dokaz Suppose $X\tau K(H,1)$ fails. This can only happen if there
is a Bockstein group $F$ in the Bockstein basis $\sigma(H)$ such
that $\dim_F(X) > 1$. That group must be either $Z_{(p)}$ or $Z/p$
for some $p$. $Z_{(p)}$ belongs to $\sigma(H)$ if and only if
$H/\Tor(H)$ is not divisible by $p$ in which case $Z_{(p)}$
belongs to $\sigma(G)$ by  \ref{ImagesBelongToClass} and
$\dim_F(X)\leq 1$ by Bockstein First Theorem. Therefore $F=Z/p$
which means that $\Tor(H)$ is not divisible by $p$. Now, either
$G$ is not divisible by $p$ or its torsion group is not divisible
by $p$ implying $\dim_F(X)\leq 1$, a contradiction.
 \hfill $\blacksquare$

\begin{Cor} \label{CohDimComparisonGvsAbG}
Let $G$ be a nilpotent group with Abelianization $\Ab(G)$ and let
$X$ be a compactum. If

\begin{itemize}
\item[a.] $X\tau K(\Ab(G),1)$,

\item[b.] $X\tau K(Q,1)$,

\item[c.] $X\tau K(Z/p^\infty,1)$ for all primes $p$,
\end{itemize}

then $X\tau K(G,1)$.
\end{Cor}
\dokaz Consider
 the lower central series of $G$: $G = \Gamma^1G\supset \Gamma^2 G\supset \ldots \supset \Gamma^i G\supset\ldots$.
 Let $F_i = \Gamma^i G/\Gamma^{i+1}G$.
 Since there is an epimorphism from $F_i \otimes \Ab(G)$ to $F_{i+1}$,
where $\Ab(G)$ is the Abelianization of $G$,  $X\tau K(F_i,1)$ for
all $i$ by \ref{EpiAndCohDim}. We proceed by induction on $c-i$
($c$ being the nilpotency class of $G$) showing that $X\tau
K(\Gamma^i G,1)$. If $c-i$ is $0$, then $\Gamma^i G=F_i$ and we
are done. Since the sequence $1 \to \Gamma^{i+1}G\to \Gamma^i G\to
F_i\to 1$ is exact, one uses a fibration $K(\Gamma^{i+1}G,1)\to
K(\Gamma^i G,1)\to K(F_i,1)$ to conclude $X\tau K(\Gamma^i G,1)$
given $X\tau K(\Gamma^{i+1} G,1)$. That constitutes the inductive
step and, as $\Gamma^1 G=G$, we get $X\tau K(G,1)$. \hfill
$\blacksquare$

\begin{Cor} \label{CohDimWRTNilpotentDp}
Let $p$ be a natural number and let ${\mathcal D}_p$ be the class
of groups $G$ such that $H/Tor(H)$ is $p$-divisible, where $H$ is
the Abelianization of $G$. If $G\in {\mathcal D}_p$ is nilpotent
and $X\tau K(Z[\frac{1}{p}]\oplus Z/p,1)$, then $X\tau K(G,1)$ for
any compactum $X$.
\end{Cor}

\begin{Cor} \label{CohDimWRTNilpotentTor}
Let $G$ be a nilpotent group. If the Abelianization $\Ab(G)$ of
$G$ is a torsion group and $X\tau K(\bigoplus\limits_{p}Z/p,1)$,
then $X\tau K(G,1)$ for any compactum $X$.
\end{Cor}

\section{Appendix B}

In this Appendix we prove results allowing us to detect homology
via maps to finite complexes with finite homotopy groups.

\begin{Lem} \label{LevinHomologyInFullStrength}
Let $A$ be a finite CW complex and $\alpha\in H_k(A)$ a nontrivial
element where $k\ge 2$. There exists a finite $(k-1)$-connected CW
complex $B$ with finite homotopy groups and a map $f\colon A\to B$
with $\beta=f_*(\alpha)$ nontrivial. Furthermore, if $\alpha$ is
of infinite order in $H_k(A)$, then $\beta$ may be assumed to be
of order $r$ for any given natural $r\ge 2$.
\end{Lem}
\dokaz With the exception of the statements about the
connectedness and the order, this is precisely Lemma 2.1 of Levin
\cite{Lev}. In the course of proving the cited lemma, Levin
constructs a ($k-1$)-connected complex $L$, and he makes $\beta$
of order $2$ if $\alpha$ has infinite order. The generalization to
arbitrary $r$ is trivial. \hfill $\blacksquare$

\begin{Lem} \label{LevinStyleHomologyLemma}
Let $M$ be a finite CW complex, and let $\PPP$ be a nonempty set
of primes. Let $\alpha\in H_k(M;Z_{(P)})$ be a nontrivial element
for some $k\ge 2$. Then there exists a finite $(k-1)$-connected CW
complex $N$ with $\PPP$-torsion homotopy groups and a map $f\colon
M\to N$ with $f_*(\alpha)$ nontrivial.
\end{Lem}
\dokaz The assumption is that there exists an element $\alpha\in
H_k(M)$ which is either $\PPP$-torsion or it has infinite order.
We can apply Lemma \ref{LevinHomologyInFullStrength} to obtain a
($k-1$)-connected finite complex with finite homotopy groups $N'$
and a map $f'\colon M\to N'$ with $\beta'=f'_*(\alpha)$ nontrivial
of order all of whose prime divisors belong to $\PPP$. Let $N'\to
N$ be localization at the set $\PPP$. Then $\beta'$ will map to
nontrivial $\beta$ under localization $\Tilde H_*(N')\to\Tilde
H_*(N)=\Tilde H_*(N')\otimes Z_{(\PPP)}$ and $N$ is (homotopy
equivalent to) the finite complex as in the statement of the
lemma. \hfill $\blacksquare$

\end{document}